\documentclass[14pt,landscape]{article}

\usepackage[colorlinks]{hyperref}
\usepackage[table]{xcolor}

\definecolor{dgreen}{rgb}{0.1, 0.4, 0.2} 
\definecolor{a}{rgb}{0.65, 0.04, 0.37}

\usepackage{lscape}
\usepackage[utf8]{inputenc}
\usepackage{graphicx}
\usepackage{color}
\usepackage{shadow}
\usepackage{geometry}
\usepackage{amsfonts,amsmath,amssymb,amsthm,nicefrac} 

\usepackage{caption}
\usepackage{subcaption}

\usepackage{mathdots}
\usepackage{amsmath}
\usepackage[table]{xcolor}
\usepackage{euscript}
\usepackage{longtable}
\usepackage{array}
\usepackage{mathtools}
\usepackage{rotating}

\begin{document}

\title{\sc
{One other parameterization of    
{$\mathrm{SU(4)}$ group} 
}}
\author{A.~Khvedelidze$^{1,2,3}$\,,  
D.~Mladenov$^{4}$\,, and A.~Torosyan$^{3,5}$
\\[1cm]
{\small $^{1}$ 
\emph{A.~Razmadze Mathematical Institute, 
Iv.~Javakhishvili Tbilisi State University, 
}}\\[-0.3cm]{\small\emph{Tbilisi, Georgia}}\\
{\small $^{2}$ 
\emph{Institute of Quantum Physics and Engineering Technologies,  Georgian Technical University,}}
\\[-0.3cm]{\small\emph{Tbilisi, Georgia}}\\
{\small $^{3}$ 
\emph{Meshcheryakov Laboratory  of Information Technologies, Joint Institute for Nuclear Research,}}
\\[-0.3cm]{\small\emph{Dubna, Russia}}\\
{\small $^{4}$ 
\emph{Faculty of Physics, 
Sofia University ``St. Kliment Ohridski'', Sofia, Bulgaria}}\\ 
{\small $^{5}$ 
\emph{A.I.~Alikhanyan National Science Laboratory (YerPhI), Yerevan, Armenia}} 
}

\date{\small \today}
\maketitle
\newpage
\begin{abstract}
We propose a special decomposition of the Lie $\mathfrak{su}(4)$ algebra into the direct sum of orthogonal subspaces,
$\mathfrak{su}(4)=\mathfrak{k}\oplus\mathfrak{a}\oplus\mathfrak{a}^\prime\oplus\mathfrak{t}\,,$ with $\mathfrak{k}=\mathfrak{su}(2)\oplus\mathfrak{su}(2)$ and a triplet of 3-dimensional Abelian subalgebras $(\mathfrak{a}, \mathfrak{a}^{\prime}, \mathfrak{t})\,,$ such that the exponential mapping of a neighbourhood of the $0\in \mathfrak{su}(4)$ into a neighbourhood of the identity of the Lie group provides the following  factorization of an element of $SU(4)$  
\[
g = k\,a\,t\,, 
\]
where 
$k \in \exp{(\mathfrak{k})} = SU(2)\times SU(2) \subset SU(4)\,,$ the diagonal matrix $t$ stands for an element from the maximal torus $T^3=\exp{(\mathfrak{t})},$ 
and the factor $a=\exp{(\mathfrak{a})}\exp{(\mathfrak{a}^\prime)}$
corresponds to a point in the double coset $SU(2)\times SU(2)\backslash SU(4)/T^3.$ 

Analyzing the uniqueness of the inverse of the above exponential mappings, we establish a logarithmic coordinate chart of the $SU(4)$ group manifold comprising 6 coordinates on the embedded manifold  $ SU(2)\times SU(2) \subset SU(4)$ 
and 9 coordinates on three copies of the regular octahedron with the edge length $2\pi\sqrt{2}\,$.  
\end{abstract}

\newpage

\tableofcontents

\newpage
%%%%%%%%%%%%%%%%%%%%%%%%%%%%%%%%%%%%%
%\listoftodos

%%%%%%%%%%%%%%%%%%%%%%%%%%%%%%%%%%%%
\section{Introduction}
%%%%%%%%%%%%%%%%%%%%%%%%%%%%%%

The paper aims to propose  one other parameterization of the special unitary group $SU(4)$ in  addition to the list  of many known ways  of describing  the unitary transformations intensively used  in pure mathematics, as well as in a diverse variety of applications~\cite{Helgason1978,Gilmore2012}.
Under ``parameterization'' we understand the term  dating back to the nineteenth century, which is  usually explicitly or implicitly  assumed in almost all relevant publications, see e.g. \cite{Murnaghan1952}-\cite{SpenglerHuberHiesmayr}.  
In our interpretation 
it means the expression of elements of a group in terms of real variables. 
Due to the non-trivial topological properties of $SU(n)$, the existence of such a unique set of parameters describing  all elements of the unitary  group is impossible. 
Therefore, it is commonly accepted that the term ``parameterization'' implies the enumeration of  ``almost all'' group elements, except those whose relative measure is zero.
\footnote{The described above approach to a group parameterization is a light version of a proper but cumbersome description of Lie groups as smooth manifolds with a certain atlas. For small values of $n$ such  chart ``coordinatization'' of a unitary group can be done explicitly, see e.g. the atlas construction for $SU(3)$ given in Appendix A of \cite{Charzynski2005}.}

Any implementation of the idea of parameterization of the operations of a Lie group is, as a rule, complemented by natural simplicity demands.
As practice shows, one's choice of coordinate charts can either lead to a drastic simplification of the whole computational scheme or ruin its efficiency at all. 
A very common example confirming the above observation is related to  situations when Lie symmetries of a physical system have to be taken into account explicitly.
 
\paragraph{Motivation of Lie group parameterizations adjusted to the description of cosets.} Formalizing the aforementioned example, 
consider two Lie subgroups $K$ and $H$  of a certain Lie symmetry group $G$ 
acting on the set of states of a physical system. As a result, the system description can be reduced to the computational problems on the equivalent classes, the double coset  $K\backslash G /H \,.$  
\footnote{Very often the double coset space is considered  for a compact Lie group
$G$ and closed connected subgroups $K$ and $H$. In this case, if  $K$ is acting freely on the coset space $G/H$, the double coset  $K\backslash G /H $ is termed as \textit{biquotient space}.
}
Hence, the symmetry urges a certain factorization of elements of $G$ adapted to its subgroup structure.
This idea turns out to be a very constructive guide for studying the generic properties of Lie groups, as well as for applying them 
\cite{Helgason1978,Gilmore2012,EdelmanJeong2022}. 
In the present article, adhering to this ideology, we consider an example of the special unitary group $SU(4)$. This choice is motivated by our recent studies of the quantum properties of a pair of two-level systems 
(2-qubits).
%%%%%%%%
\footnote{Shortly, when describing 2-qubits, the double coset $SU(2)\times SU(2)\backslash SU(4)/T^3 $ arises as follows. From one side, 
the coset $SU(4)/T^3$ appears since a 2-qubit state as a whole  is characterized by  $SU(4)$ invariants,  the spectrum of  density matrix and diagonalizing unitary transformation, up to isotropy group of state.  The latter for a generic case of non-degenerate  spectrum of density matrices is $T^3\,.$
From the other side,
the subgroup 
$SU(2)\times SU(2)\subset SU(4)$ appears as the subgroup characterizing individual qubits as elementary objects. }
%%%%%%%
Further on, we introduce the parameterization of $SU(4)$ adjusted to the double coset corresponding to its subgroups, 
$K= SU(2)\times SU(2)$ and maximal torus $H=T^3$.
More precisely, in the next section we describe the following factorization  of  element $g \in SU(4)$\,,  abbreviated  afterwards  as $K\mathcal{A}T\--$decomposition, 
\begin{equation}
\label{eq:SU(4)coord}
g= k\,a\,t \,,
\end{equation}
where $k \in K, \, t \in T^3 $ and factor 
$a$ corresponds to the double coset $SU(2)\times SU(2)\backslash SU(4)/T^3\,.$ 
The factorization 
(\ref{eq:SU(4)coord}) is an addition to the family of different parameterizations of compact Lie groups well adapted to the descriptions of their cosets and double cosets.
Here we only point out some of them, the famous KAK decomposition of a Lie group based on the single Cartan involution \cite{Helgason1978,Gilmore2012,EdelmanJeong2022}, pairs of involutions \cite{Matsuki1995,Matsuki1997,Miebach2007} or Kobayashi triples  $(U(n),\,  U(p)\times U(q),\,
U(n_1)\times U(n_2)\dots U(n_l))\,$ \cite{Kobayashi2007}.
A detailed classification of existing factorization schemes of unitary groups can be found in the review article  \cite{EdelmanJeong2023}.

Our paper  is  organized as follows. 
We start with a  brief introduction of the idea of canonical coordinates on the Lie group manifold which is based on the fundamental result regarding the exponential mapping from the Lie algebra to the corresponding Lie group. In order to introduce nonstandard coordinates on $SU(4)$, we describe the adjusted decomposition of the algebra $\mathfrak{su}(4)$ into the direct sum of orthogonal vector spaces constructed as a span of five non-overlapping subsets of the tensorial Fano  basis of $\mathfrak{su}(4)$\,.
Performing the exponential mapping of the corresponding components of decomposition, we present the explicit form of  all three factors in $K\mathcal{A}T$\--decomposition (\ref{eq:SU(4)coord}). 
In Section \ref{sec:LogChart}, the analyticity domain of the introduced parameterization will be found. The relation of the parameters to the angles of rotation in the six-dimensional Euclidean space is established considering the group $SU(4)$ as double cover of the special orthogonal group $SO(6)\,.$ 
The article is concluded by a discussion on factorising the elements of $SU(4)$ into elements of maximal Abelian subgroups and special Hadamard matrices of the 4-th order. 
We also comment on the ways to use the proposed parameterization of $SU(4)$ to describe the entanglement of 2-qubit states.

%%%%%%%%%%%%

\section{Method of canonical coordinates for the $SU(4)$ group}
%%%%%%%%%%%%%%
%%%%%%%%%%%%%%%%%%%%""""""""""""""""""""

\paragraph{Canonical charts on Lie groups.}
Various forms of parameterizations of the Lie group manifold ensue from the structure of the associated
Lie algebra. 
One of the basic theorems of the Lie group theory (see e.g. 
\cite{Postnikov1986}, Lecture 4, 
\cite{Lee2013}, Theorem 20-12 and Problem 20-3)  asserts: 

\subparagraph{Theorem on canonical charts.}
\textit{If the Lie algebra
$\mathfrak{g}$ of the Lie group $G$ is decomposed into the direct sum of any number of subspaces, 
\begin{equation}
\label{eq:decalgebra}
\mathfrak{g}=\mathcal{A}_1\oplus\mathcal{A}_2
\oplus \dots \oplus \mathcal{A}_m\,,
\end{equation}
then the mapping $ \Phi: \mathfrak{g} \to G,$
\begin{equation}
\label{eq:mapping}
   \Phi(X)
   = \exp a_1 \cdot \exp a_2 {\cdot} 
\dots {\cdot} \exp a_m\,,
\end{equation} 
is a diffeomorphism from some neighbourhood of $(0,0, \dots, 0)$ in 
$\mathfrak{g}$ to a neighbourhood of the identity element of group $G$\,, 
assuming that $a_i \in \mathcal{A}_i$\,, $i=1,2,\dots m $\,, are components of vector $X$ in the decomposition (\ref{eq:decalgebra}) for any $X \in\mathfrak{g}$\,. 
}

Relation  (\ref{eq:mapping}) defines  certain types of charts called  
\textit{canonical charts}  with corresponding  \textit{canonical coordinates}
on group manifold $G$. Particularly  
for $m=1$, the  decomposition (\ref{eq:decalgebra}) defines the \textit{canonical coordinates of the first kind}, while  for the other extreme
case, $m=\dim T_e(G)$, i.e.,  when all subspaces $\mathcal{A}_i$ are one-dimensional, the canonical chart  is parameterized by the \textit{canonical coordinates of the second  kind}.

Following this basic observation, bearing in mind our aim to describe a coordinate system on the group manifold adapted to the left and right action of the $\mathrm{SU(2)}\times\mathrm{SU(2)}$ and $\mathrm{T}^3$ subgroups,  we will use  
exponential mapping providing the sought-for parameterization  of $SU(4)$. The later originates in the process of eliciting special subspaces in the Lie algebra $\mathfrak{su}(4)$ \---  the subalgebra $\mathfrak{su}(2)\oplus \mathfrak{su}(2)$ and three copies of 3-dimensional Abelian subalgebras. After exponential mapping the complemented subspace affords the introduction of a canonical chart with coordinates  on the double coset $\mathrm{SU(2)\times SU(2)}\backslash \mathrm{SU(4)}/\mathrm{T}^3$\,.

%%%%%%%%%%%%%
\subsection{KAT-decomposition of the Lie algebra $\mathfrak{su}(4)$}

As a first step of the outlined program, we shall prove the following statement.

\noindent{\bf Proposition I:\,} The $\mathfrak{su}(4)$ algebra admits of decomposition into the following  direct sum of subspaces:
\begin{equation}
\label{eq:NewDecomp}
  \mathfrak{su}(4) = \mathfrak{k}
  \oplus  \mathfrak{a}\oplus
\mathfrak{a}^\prime\oplus \mathfrak{t}\,,
\end{equation}
where $\mathfrak{t}$ is Cartan subalgebra of $\mathfrak{su}(4)$, $\mathfrak{k}:=\mathfrak{su}(2)\oplus \mathfrak{su}(2)\,,$ 
and components  $\mathfrak{a}$ and $\mathfrak{a}^\prime$ in (\ref{eq:NewDecomp}) are 3-dimensional Abelian subalgebras, such that 
\begin{equation}
\label{eq:com1}
[\mathfrak{a}^\prime\,, \mathfrak{a}] \subseteq \mathfrak{k}\,.
\end{equation}
Furthermore, the elements of decomposition  in 
(\ref{eq:NewDecomp})
obey the following commutation relations:
\begin{eqnarray}
\label{eq:com2}
&&[\mathfrak{k}, \mathfrak{k}] \subseteq \mathfrak{k}\,, \qquad 
[\mathfrak{t}, \mathfrak{a} ] \subseteq \mathfrak{k}\,,
\qquad 
[\mathfrak{t}, \mathfrak{a}^\prime ] \subseteq \mathfrak{k}\,,
\qquad
[\mathfrak{k}, \mathfrak{t}] \subseteq \mathfrak{a}\oplus
\mathfrak{a}^\prime\,.
\end{eqnarray}
%
%%%%%%%%%%%%%%%%%%%%%%%%%%%

To prove this  proposition,  we construct  every component of the direct sum 
(\ref{eq:NewDecomp}) as a span of orthogonal subsets of a special basis of the $\mathfrak{su}(4)$\ algebra. 
With this goal, we use  the Fano  basis of the $\mathfrak{su}(4)$ algebra constructed via tensor products of the form 
\[
\sigma_{\mu\nu}=\frac{1}{{2\imath}}\,\sigma_{\mu}\otimes \sigma_\nu\,,
\]
where  4-tuple $\sigma_\mu=(\mathbb{I}_2\,, \boldsymbol{\sigma})$ consists from  the unit $2\times 2$ matrix $\mathbb{I}_2$ and the Pauli matrices
$\boldsymbol{\sigma}=(\sigma_1, \sigma_2,
\sigma_3)$. 
\footnote{
The elements of the Fano basis are explicitly given in the Appendix.
}
Afterwards, for convenience,  a single index notation  for 
$\mathfrak{su}(4)$ algebra basis  $\boldsymbol{\lambda}=\{\lambda_1, \dots, \lambda_{15} \}$ will be used, assuming the following order of basis elements: 
\begin{equation}
\label{eq:basisLambdaF}
\boldsymbol{\lambda} = 
\{\sigma_{10},
\sigma_{20}, \sigma_{30},
\sigma_{01}, \sigma_{02}, \sigma_{03}, 
\sigma_{11},
\sigma_{12}, \sigma_{13}, \sigma_{21}, \sigma_{22}, \sigma_{23},
\sigma_{31}, \sigma_{32}, \sigma_{33}\}\,.
\end{equation}

At first, we construct the direct sum decomposition 
of Lie algebra  $\mathfrak{su}(4)$ considered as a real vector space and then  check  the corresponding commutation relations.
In order to find the decomposition  (\ref{eq:NewDecomp}), let us split the basis $\boldsymbol{\lambda}$ into five complementary  subsets:
\[\boldsymbol{\lambda}=
\bigcup_{i=1}^5\boldsymbol{\Lambda}_i\,,
\]
where 
$
\boldsymbol{\Lambda}_1=\{\lambda_1, \lambda_{4}, \lambda_7\}, \quad 
\boldsymbol{\Lambda}_2 =\{\lambda_9, \lambda_{11}, \lambda_{13}\}, \quad
\boldsymbol{\Lambda}_3=\{\lambda_2, \lambda_{8}, \lambda_{14}\}, \quad 
\boldsymbol{\Lambda}_4 =\{\lambda_5, \lambda_{10}, \lambda_{12}\}$ and $ 
\boldsymbol{\Lambda}_5=\{\lambda_3, \lambda_{6}, \lambda_{15}\}\,.
$
Defining  the corresponding spans, $\mathcal{A}_i=\mbox{span}_{\mathbb{R}}(\boldsymbol{\Lambda}_i), \,  i=1,2,\dots, 5\,, $ we see that the  direct summation  of spans $\mathcal{A}_i$ gives the algebra  
\[
\mathfrak{su}(4)=\bigoplus_{i=1}^5\, \mathcal{A}_i  
\]
due to the specific form of sub-basis tuples 
$\boldsymbol{\Lambda}$\,.
Finally, we can directly verify that  relations (\ref{eq:com1}) and (\ref{eq:com2})
hold
(c.f. Table \ref{T:CommRelsu(4)}. in Appendix)  if the elements in the sought-for decomposition (\ref{eq:NewDecomp})
are identified with the above spans as:
\[
\mathfrak{a}=\mathcal{A}_1\,, \qquad
\mathfrak{a}^\prime=\mathcal{A}_2\,,\qquad
\mathfrak{k}=\mathcal{A}_3\oplus \mathcal{A}_4 
\,,
\qquad 
\mathfrak{t}=\mathcal{A}_5\,. 
\]

%%%%%%%%%
\subsection{KAT-parameterization of the $SU(4)$ group}

Applying the theorem on the canonical chart for the Lie group (\ref{eq:mapping}) and the direct sum decomposition (\ref{eq:NewDecomp}), we arrive at our main  Proposition. 

\noindent{\bf Proposition II:\,}
\textit{The decomposition  of the algebra  $\mathfrak{su}(4)$ into the direct sum defined in Proposition I determines, via  exponential map  $\exp : \mathfrak{su}(4) \to \mathrm{SU(4)}$\,,
the following factorization form of the $SU(4)$ group element:
\begin{equation}
\label{eq:SU4groupKAA'T}
g= \exp{(\mathfrak{k})}\cdot\exp{(\mathfrak{a})}\cdot
\exp(\mathfrak{a}^\prime)\cdot\exp{(\mathfrak{t})}\,, \qquad \ g\in SU(4)\,.
\end{equation}
In (\ref{eq:SU4groupKAA'T})
$exp{(\mathfrak{k})}$  corresponds to the subgroup $  
K=\exp(\mathfrak{su}(2))\times
\exp(\mathfrak{su}(2))\, \subset \mathrm{ SU(4)}$, while 
$\exp{(\mathfrak{t})}$  is the maximal torus $T^3$ in $\mathrm{SU(4)}$.
}

Coefficients of the expansion of $\mathfrak{k}, \mathfrak{a}, \mathfrak{a}^\prime $ and $\mathfrak{t}$
over the basis of Lie algebra $\mathrm{su}(4)$ are identified  with a point in $\mathbb{R}^{15}$ and, therefore, mapping (\ref{eq:SU4groupKAA'T}) defines a canonical chart on group manifold $SU(4)$ in the vicinity of identity. 
However, since every Lie group is a parallelizable manifold, the tangent space at the identity can be translated over the group, which allows to extend the chart to ``almost'' the whole group.
To formalize this observation, we write the following K$\mathcal{A}$T-decomposition for the $SU(4)$ group: 
\begin{equation}
\label{eq:SU4KAT}
    SU(4) = K\mathcal{A}T\,.  
\end{equation}
Below, in order to justify the correctness of the suggested  parameterization (\ref{eq:SU4KAT}),  we will analyze in detail the domain of definition of the canonical coordinates. With this aim, let us begin with explicit representations of each of the three factors in the $K\mathcal{A}T$\--decomposition
of $SU(4)$ in the introduced  canonical chart.  

%%%%%%%%%%%%%%%%%%
\subsection{Explicit formulae for KAT\--decomposition factors}
%%%%%%%%%%%%%%%

\subparagraph{
K-factor.}
In the proposed $K\mathcal{A}T$-decomposition of the $SU(4)$ group, the $SU(2)\times SU(2)$ factor  is given through the exponential map of $\mathfrak{k}=\mathrm{span}_{\mathbb{R}}(\lambda_2, \lambda_8, \lambda_{14}, \lambda_5, \lambda_{10}, \lambda_{12})\,,$
\begin{equation}
   \label{eq:K1}
K=\exp(\mathfrak{k})=\exp\left((\boldsymbol{u},\boldsymbol{\Lambda_3}) + (\boldsymbol{v},\boldsymbol{\Lambda_4})\right)\,,
\end{equation}
with real 3-tuples, 
$\boldsymbol{u}=\{u_2, u_8, u_{14}\}$
and $\boldsymbol{v}=\{v_5, v_{10}, v_{12}\}\,.$

On the other hand, 
the standard embedding of $\mathfrak{su}(2)\oplus\mathfrak{su}(2)$ is realized as  the composition of two types of 
$\mathfrak{su}(2) \hookrightarrow \mathfrak{su}(4)$ embeddings:
\begin{equation}
\label{eq:Embed1}
    (I):\qquad \qquad
\mathfrak{su}(2) \hookrightarrow \mathfrak{su}(4): \quad \left[
\begin{array}{ccc}
x        & z \\
\bar{z}  & -x
\end{array}
\right]
\hookrightarrow 
\left[
\begin{array}{cc|cc}
x       & ~~0~        & z  & 0  \\
0       & ~~x~ \       & 0  & z  \\ \cline{1-4}
\bar{z} &  ~~0        & -x &  0 \\
0       & ~~\bar{z}   &   0 & -x
\end{array}
\right]\,,
\end{equation}
and 
\begin{equation}
\label{eq:Embed2}
(II):\qquad \qquad
\mathfrak{su}(2) \hookrightarrow \mathfrak{su}(4): \quad \left[
\begin{array}{ccc}
y        & w\\
\bar{w}  & -y
\end{array}
\right]
\hookrightarrow 
\left[
\begin{array}{cc|cc}
y       & w   & 0  & 0  \\
\bar{w} & -y  & 0  & 0  \\ \cline{1-4}
 0      &  0  & y &  w \\
0       & 0   &   \bar{w} & -y
\end{array}
\right]\,.
\end{equation}
In terms of the Fano basis, 
this embedding results in a subalgebra  different from $\mathfrak{k}$, namely Lie(SU(2)$\times $SU(2)) =
$\mathrm{span}_{\mathbb{R}}(\lambda_1, \lambda_2, \lambda_{3}, \lambda_4, \lambda_{5}, \lambda_{6})\,.$
The exponentiation of this subalgebra  gives a conventional embedding of $SU(2)\times SU(2)$ into $SU(4)$:
\begin{equation}
\label{eq:convenSU2XSU2}
U(\varphi, \boldsymbol{n})\otimes U(\psi, \boldsymbol{m})=
\exp\left(-\frac{1}{2}\varphi\,(\boldsymbol{n}\cdot\boldsymbol{\sigma})\right)
\otimes
\exp\left(-\frac{1}{2}\psi(\boldsymbol{m}\cdot\boldsymbol{\sigma})\right) \subset 
SU(4)
\,,
\end{equation}
describing two independent rotations on angles $\varphi$ and $\psi$ around two unit axes $\boldsymbol{n}$
and $\boldsymbol{m}$, respectively.

In order to connect  $K\--$factor  (\ref{eq:K1}) with the conventional 
embedding   (\ref{eq:convenSU2XSU2}),
we note the following relations between the elements of the Fano basis:  
\begin{equation}
\label{eq:Rtr}
    \mathrm{R}^\dagger\{ \lambda_1, \lambda_2, \lambda_3,
    \lambda_4, \lambda_{5}, \lambda_{6}\}
    \mathrm{R} = \{ \lambda_5, \lambda_{12}, \lambda_{10}, 
    \lambda_{14},  - \lambda_{2}, \lambda_{8}\}\,,
\end{equation}
 where $R$ denotes the ``magic'' $4\times 4$ unitary matrix \cite{FujiiSuzuki2007}: 
\begin{equation}
    \mathrm{R}=\frac{1}{\sqrt{2}}\,\left(
\begin{array}{cccc}
 1 & 0 & 0 & -{i} \\
 0 & -{i} & -{1} & 0 \\
 0 & -{i} & {1} & 0 \\
 {1}& 0 & 0 & {i} \\
\end{array}
\right)\,. 
\end{equation}
Using (\ref{eq:Rtr}), one can easily  find expressions  of conventional coordinates,  rotation angles  
$\varphi, \psi$ and rotation axes $(\boldsymbol{n}, \boldsymbol{m} )$ of  two copies of the SU(2) group
in terms of the coordinates $\boldsymbol{u}$ and $\boldsymbol{v}$:
\begin{equation}
    \label{eq:uvphipsi}
    \varphi\,\boldsymbol{n}= 
    \{v_5, v_{10}, v_{12}\}\,,\qquad 
 \psi\,\boldsymbol{m} = 
    \{-u_2, u_{8}, u_{14}\}\,.   
\end{equation}
Hence, we  arrive at the representation for the $K\--$factor: 
\begin{equation}
\label{eq:K2}
K =\mathrm{R}^\dagger\left(
\exp(-\frac{1}{2}\varphi\,(\boldsymbol{n}\cdot\boldsymbol{\sigma}))\otimes  
\exp(-\frac{1}{2}\psi\,(\boldsymbol{m}\cdot\boldsymbol{\sigma}))
\right)
\mathrm{R}\,. 
\end{equation}

\subparagraph{$\mathcal{A}$-factor.}
Let us identify $\mathfrak{a}$ and $\mathfrak{a}^\prime$ with  $\mathbb{R}^6$ by introducing two triplets of coordinates 
$\boldsymbol{\alpha}=(\alpha_1, \alpha_2,\alpha_3)$ and $\boldsymbol{\beta}=(\beta_1, \beta_2,\beta_3)$ via the following decompositions:  
\[
(\boldsymbol{\alpha}, \boldsymbol{\Lambda}_1)=
\alpha_1\lambda_1+\alpha_2\lambda_4+\alpha_3\lambda_7
\in \mathfrak{a} 
\qquad\mbox{and}\qquad
(\boldsymbol{\beta}, \boldsymbol{\Lambda}_2)
=\beta_1\lambda_9+\beta_2\lambda_{11}+
\beta_3\lambda_{13}
\in
\mathfrak{a}^\prime\,.
\]
The corresponding $\mathcal{A}$-factor in the 
$K\mathcal{A}T$\--decomposition represents the product of  two exponents:   
\begin{equation}
\mathcal{A}_1= 
\exp{\displaystyle{( \boldsymbol{\alpha}, \boldsymbol{\Lambda}_1)}}\,,
\qquad 
\mathcal{A}_2 =
\exp{( \boldsymbol{\beta}, \boldsymbol{\Lambda}_2)}\,.
\end{equation}
Each of the exponents expands over the corresponding elements of the Fano basis: 
\begin{eqnarray}
\label{eq:M1M21}
\mathcal{A}_1&=&x_0\mathbb{I}_{4}+ x_1\lambda_1+
x_2\lambda_4+x_3\lambda_7\,,\\
\mathcal{A}_2&=&y_0\mathbb{I}_{4}+
  y_1\lambda_9+
y_2\lambda_{11}+y_3\lambda_{13}\,. 
\label{eq:M1M22}
\end{eqnarray}
The coefficients $x_0, x_1, x_2$ and $ x_3$  in the expansion of the factor $\mathcal{A}_1$ are complex valued functions of the angles $\alpha_1\,, \alpha_2$ and $\alpha_3$:
\begin{eqnarray}
  x_0&= &2\imath\left(\cos \left(\frac{\alpha _1}{2}\right) \cos \left(\frac{\alpha _2}{2}\right) \cos \left(\frac{\alpha _3}{2}\right)+i \sin \left(\frac{\alpha _1}{2}\right) \sin \left(\frac{\alpha _2}{2}\right) \sin \left(\frac{\alpha _3}{2}\right)\right)\,,\\
  x_1&=&2\left(\sin \left(\frac{\alpha _1}{2}\right) \cos \left(\frac{\alpha _2}{2}\right) \cos \left(\frac{\alpha _3}{2}\right)-i \cos \left(\frac{\alpha _1}{2}\right) \sin \left(\frac{\alpha _2}{2}\right) \sin \left(\frac{\alpha _3}{2}\right) \right)\,,\\
  x_2&=&2\left(\cos \left(\frac{\alpha _1}{2}\right) \sin \left(\frac{\alpha _2}{2}\right) \cos \left(\frac{\alpha _3}{2}\right)-i \sin \left(\frac{\alpha _1}{2}\right) \cos \left(\frac{\alpha _2}{2}\right) \sin \left(\frac{\alpha _3}{2}\right) \right)\,,
  \\
  x_3&=& 2\left(\cos \left(\frac{\alpha _1}{2}\right) \cos \left(\frac{\alpha _2}{2}\right) \sin\left(\frac{\alpha _3}{2}\right) - i \sin \left(\frac{\alpha _1}{2}\right) \sin \left(\frac{\alpha _2}{2}\right) \cos \left(\frac{\alpha _3}{2}\right)\right)\,,
\end{eqnarray}
while coefficients 
$y_0, y_1, y_2$ and $ y_3$ of the expansion of the factor $\mathcal{A}_2 $ are, respectively, the same functions of $\boldsymbol{\beta}\,.$
As a result, factor  $\mathcal{A}$ is decomposing as 
\begin{equation}
\label{eq:decomM}
\mathcal{A} =  \sum_{A=0}^{15}\,M_A\lambda_A\,,
\end{equation}
with the following components of the 16-tuple:
\begin{eqnarray}
\label{eq:16vector}
M_{0}&= &-\frac{\imath}{2}\,x_0y_0\,,\\ 
M_i&=&-\frac{\imath}{2}\left(
x_1y_0\,,-\imath x_3y_3\,,  x_2 y_3\,,
x_2y_0\,, -\imath x_3y_1\,,  x_1y_1\,,
\right. \\
&&
\left.
x_3y_0\,, -\imath x_2y_1\,, x_0y_1,
-\imath x_1y_3\,, x_0y_2\,, \imath x_2y_2\,, 
x_0y_3\,, \imath x_1y_2\,, -x_3y_2
\right)\,.
\end{eqnarray}
It is worth noting that Abelian factors $\mathcal{A}_1$ and $\mathcal{A}_2$ in the $K\mathcal{A}T$\--decomposition can be diagonalized by the corresponding orthogonal transformations: 
\begin{eqnarray}
\label{eq:Afactor1}
 \mathcal{A}_1 &=& \exp{\displaystyle{( \boldsymbol{\alpha}, \boldsymbol{\Lambda}_1)}} = H_1 \,D(\boldsymbol{\alpha})\, H_1^T\,,
 \\
 \mathcal{A}_2 &=& \exp{( \boldsymbol{\beta}, \boldsymbol{\Lambda}_2)} = H_2\,
 D(\boldsymbol{\beta}) \,H_2^T\,, 
 \label{eq:Afactor2}
\end{eqnarray}
where 
\begin{equation}
\label{eq:S}
H_1=
\frac{1}{2} \, \left(
\begin{array}{cccc}
 1 & 1 & -1 & -1 \\
 1 & -1 & -1 & 1 \\
 1 & -1 & 1 & -1 \\
 1 & 1 & 1 & 1 \\
\end{array}
\right)\,,
\qquad
\qquad 
H_2=\frac{1}{2} \, \left(
\begin{array}{cccc}
 -1 & 1 & -1 & 1 \\
 -1 & 1 & 1 & -1 \\
 -1 & -1 & 1 & 1 \\
 1 & 1 & 1 & 1 \\
\end{array}
\right)\,,
\end{equation}
and their diagonal counterparts are given by means of the matrix function $D(\boldsymbol{x})$ of 3-tuple $\boldsymbol{x}$: 
\begin{equation}
\label{eq:diag}
D(\boldsymbol{x})=\mbox{diag}
||\,
e^{-\frac{\imath}{2}\left(x_1+x_2+x_3\right)}\,, 
e^{\frac{\imath}{2} \left(x_1+x_2-x_3\right)}\,, 
e^{\frac{\imath}{2} \left(x_1-x_2+x_3\right)}\,, 
e^{\frac{\imath}{2}\left(-x_1+x_2+x_3\right)}\,
||\,.
\end{equation}
Moreover,  both matrices ${H}_1$ and $H_2$
are Hadamard  $4\times 4$ matrices.
They are 
equivalent to the standard Hadamard matrix, constructed using Sylvester's method from the $2\times2$ Hadamard matrix as follows:   
\begin{equation}
    H^{(4)}= \begin{pmatrix}
        H^{(2)}& H^{(2)}\\
        H^{(2)} & - H^{(2)}
    \end{pmatrix}\,, \qquad 
H^{(2)}=\begin{pmatrix}
    1&1\\
    1&-1
\end{pmatrix}\,. 
\end{equation}
Indeed, one can show that the matrices ${H}_1$ and ${H}_2$ can be obtained from $H_4$ by a signed permutation of rows and columns, i.e. 
\begin{equation}
    H_i = \frac{1}{2}\,P_i H^{(4)} Q_i \,,\qquad 
    P_i, Q_i \in S_4\times \mathbb{Z}_2\,,
\end{equation}
where 
\begin{equation}
    P_1 = \left(
\begin{array}{cccc}
 1 & 0 & 0 & 0 \\
 0 & 1 & 0 & 0 \\
 0 & 0 & 0 & 1 \\
 0 & 0 & 1 & 0 \\
\end{array}
\right)\,, \quad 
    Q_1 = \left(
\begin{array}{cccc}
 1 & 0 & 0 & 0 \\
 0 & 1 & 0 & 0 \\
 0 & 0 & -1 & 0 \\
 0 & 0 & 0 & -1 \\
\end{array}
\right)\,, \quad 
    P_2 = \left(
\begin{array}{cccc}
 0 & 0 & 1 & 0 \\
 0 & 0 & 0 & 1 \\
 1 & 0 & 0 & 0 \\
 0 & -1 & 0 & 0 \\
\end{array}
\right)\,, \quad 
    Q_2 = \left(
\begin{array}{cccc}
 -1 & 0 & 0 & 0 \\
 0 & 0 & 0 & 1 \\
 0 & -1 & 0 & 0 \\
 0 & 0 & 1 & 0 \\
\end{array}
\right)\,. 
\end{equation}

\subparagraph{T-factor.}
Finally, exponential mapping of the last component of the direct sum decomposition (\ref{eq:NewDecomp}), the maximal Abelian algebra, 
$\mathfrak{t}=\theta_3\lambda_3+\theta_6\lambda_6+\theta_{15}\lambda_{15}\,,$
results in  3\--dimensional torus $ T^3 \subset SU(4)$ 
given by means of the diagonal matrix $D(\boldsymbol{\theta})$\,, $\boldsymbol{\theta}=\{\theta_3, \theta_6, \theta_{15}\}\,.$

The next subsection is devoted to the discussion  of  the range of parameters in  each of  those factors,  
$K, \mathcal{A}$, and $T$,  ensuring  the correctness of the description of the $SU(4)$ transformations within the given chart.

%%%%%%%%%%%%%%%""""""""""""""""""""""
\subsection{The range of coordinates in the logarithmic chart}
%%%%%%%%%%%%%%"""""""""""""""""
\label{sec:LogChart}

\paragraph{Normal neighbourhood.} According to the theorem of canonical charts, the introduced set of 15 coordinates, including 6 coordinates $(\varphi, \boldsymbol{n})$ and $(\psi, \boldsymbol{m})$ describing the $SU(2)\times SU(2)$ subgroup, 3 torus angles  $\boldsymbol{\theta}$ and 6 angles  $(\boldsymbol{\alpha}, \boldsymbol{\beta})$ describing the neighbourhood of (0) in $\mathfrak{su}(4)$ algebra, provides the parameterization of elements in the neighbourhood of the identity element of the $SU(4)$ group. 
The admissible range  of these coordinates  is dictated by the requirement of the diffeomorphism mapping that includes 
1-to-1 correspondence  between elements of the $SU(4)$ group and the moduli space of parameters. 
Following these requirements, one can find  the analyticity domain and thereby correctly define a specific  chart on the $SU(4)$ group manifold. Being related to the inverse of the exponential mapping, this chart is sometimes called \textit{logarithmic chart}.
The theoretical background justifying the parameterization of Lie groups based on identifying its Lie algebra with the corresponding Euclidean space and subsequent exponential mapping relies on the theorem determining those domains of coordinates where the restriction of the matrix exponent has a unique inverse (cf. Theorem 11.3  in \cite{Higham2008}).

\textbf{Theorem on exponential mapping.}\textit{ Let G be a compact and connected matrix Lie group. The exponential function $\exp_A $ is analytic, with an analytic inverse on a bounded open neighbourhood of the origin given by
\begin{equation}
\label{eq:anDomain}
 U=\{ A \in \mathfrak{g}\, |\, |\operatorname{Im}[\mathrm{Eigenvalue}(A)]|<\pi \}\,.  \end{equation}
}

This theorem indicates that, due to the multi-valued nature of the logarithm, it is not generally the case that $\log(\exp A)=A$.
The constraints (\ref{eq:anDomain}) generalize the result
known from the theory of complex variables on the relationship between 
$\log(e^z)$ and $z \in \mathbb{C}$, showing  its uniqueness,
$\log(e^z)=z$\,, for domain $\operatorname{Im}(z) \in (-\pi, \pi]$.

\begin{figure}[h!]
\center{
\includegraphics[scale=0.65]{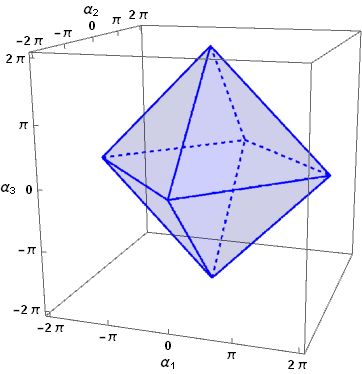}
}
\caption{One of the 3 regular octahedrons as a domain of coordinates in the logarithmic chart.  }
\label{fig:Octahedron}
\end{figure}

Using the equations (\ref{eq:anDomain}),
we are able to find  the domain  of parameters in 
(\ref{eq:SU4groupKAA'T})
corresponding to the logarithmic chart  covering almost the whole $SU(4)$. 
\footnote{
Note that for the normal coordinates the closure of the normal neighbourhood is equal to the whole group.
}
All factors of the $K\mathcal{A}T\--$decompositions display  
$4\pi$\--periodicity, inherited from the 
periodicity  of  $SU(2)$ factors and 
the Abelian subgroups,  
\begin{eqnarray}
&&U(\varphi+4\pi, \boldsymbol{n})=U(\varphi, \boldsymbol{n})\,, \quad
U(\psi+4\pi, \boldsymbol{m})=
U(\psi, \boldsymbol{m})\,,\\
&&
D(\boldsymbol{x+4\pi})= D(\boldsymbol{x})\,.
\end{eqnarray}
In order to have 1-to-1 correspondence  between 
group elements and angles, we restrict the range of all angles by the interval $(-2\pi,2\pi ]$\,, 
and hence the group elements are in correspondence with the real 11-dimensional  torus.  However,  not for all points of this torus the  constraints  (\ref{eq:anDomain})  are satisfied automatically. Thus, some of the angles need to be constrained in order to be  proper coordinates on the $SU(4)$ group manifold.
Requirements (\ref{eq:anDomain}) introduce further restrictions on the angles' domains of definition.  
Indeed, applying (\ref{eq:anDomain})  to
the factor corresponding to the subgroup
$SU(2) \times SU(2)$, we find the domain of analyticity for the
3-tuples of angles $\phi$ and $\boldsymbol{\psi}$\,:
\begin{equation}
    \label{eq:anpsiphi}
 0 \leq \varphi < 4\pi\,, \qquad  0\leq \psi < 4\pi \,,  
\end{equation}
while for the $\mathcal{A}\--$factor,
written as the product of two Abelian subgroups  
(\ref{eq:Afactor1}),  (\ref{eq:Afactor2}),  and the diagonal torus  $T^3$, the analyticity domain
for 3-tuples $\boldsymbol{\alpha}, \boldsymbol{\beta}$, and $\boldsymbol{\theta}$ reads:  
\begin{eqnarray}
&&  |\alpha_1+\alpha_2+\alpha_3| < 2\pi\,, \qquad 
 |\beta_1+\beta_2+\beta_3|< 2\pi\,,\qquad 
 |\theta_3+\theta_6+\theta_{15}|< 2\pi\,, 
 \\ 
  &&  |\alpha_1+\alpha_2-\alpha_3|< 2\pi\,, \qquad 
 |\beta_1+\beta_2-\beta_3|< 2\pi\,,
 \qquad 
 |\theta_3+\theta_6-\theta_{15}|< 2\pi\,,
 \\  
&&  |\alpha_1-\alpha_2+\alpha_3|< 2\pi\,,
\qquad 
 |\beta_1-\beta_2+\beta_3|< 2\pi\,, \qquad
 |\theta_3-\theta_6+\theta_{15}|< 2\pi\,,
   \\
&&|\alpha_2+\alpha_3-\alpha_1|< 2\pi\,,
\qquad 
|\beta_2+\beta_3-\beta_1|< 2\pi\,,
\qquad
|\theta_{3}+\theta_{15}-\theta_{3}|< 2\pi\,.
\end{eqnarray}
These inequalities describe the direct product of three 
regular octahedrons (Fig. \ref{fig:Octahedron}) whose points have the Cartesian coordinates 
$\boldsymbol{\alpha}, \boldsymbol{\beta}$, and 
$\boldsymbol{\theta}$, respectively.
The points of these octahedrons are in 1-to-1 correspondence with the group elements, and thus 
$\boldsymbol{\alpha}, \boldsymbol{\beta}$, and 
$\boldsymbol{\theta}$ together with $(\varphi, \boldsymbol{n})$ and  $(\psi, \boldsymbol{m})$ are  coordinates in the logarithmic chart defined by the $K\mathcal{A}T$\--decomposition of $SU(4)$.

%%%%%%%%%%%%%%%%%%%%%%%%%%%%%%%  End figure

\paragraph{Parameterizing $SU(4)$ and $SO(6)$.} Establishing the domain of definition of group parameters, it is worth to  take into account the existence of non-equivalent connected Lie groups whose Lie algebra is $\mathfrak{su}(4)$.
Indeed, our construction of the $SU(4)$ group chart relies  on the exponential mapping from the algebra,
and thus, different domains of definition of parameters will result in different Lie groups. 
To verify that our parameterization indeed gives elements of the $SU(4)$ group, let us recall the relation between the $SU(4)$ group and special orthogonal group $SO(6)\,,$ whose Lie algebras are isomorphic. 
It is well-known that the universal cover of  
$SO(n)$ is a simply-connected compact Lie group  $\mathrm{Spin}(n)$. Moreover, since fundamental group $\pi_1(SO(n)) = \mathbb{Z}_2$, the universal covering group of $SO(n)$ is a double covering group $\mathrm{Spin}(n)$:
\[
\mathbb{Z}_2 \rightarrow \mathrm{Spin}(n) \rightarrow SO(n)\,.
\]
This observation together with the following isomorphism:  
 \[
 \mathrm{Spin}(3)\simeq SU(2),\qquad  
\mathrm{Spin}(4) \simeq  SU(2) \times SU(2),\qquad \mbox{and}\qquad   \mathrm{Spin}(6) \simeq  SU(4)\,,
\]
allows us to identify subsets of group parameters in the proposed $K\mathcal{A}T$\--factorization  with the subsets of angles describing the $SO(6)\,$ rotations.  For this reason, let us write down an explicit formulae for the mapping between $SO(6)$ and $SU(4)$ based on the isomorphism between the Lie algebras $\mathfrak{su}(4)$ and $\mathfrak{so}(6)\,.$

\paragraph{Isomorphism between $SU(4)$ and $\mathrm{Spin}(6)$\,.}
Let $L_{(i,j)}$ be  15 standard  generators  of the Lie algebra of the Lie group $SO(6)$, i.e., the 
skew-symmetric $6\times 6$ matrices, each with a single 1 entry above the diagonal and a corresponding -1  entry below the diagonal (all other entries are zero):
\begin{equation}
\label{eq:basisso(6)}
\left(L_{(i,j)}\right)^{mn}=\delta^m_i\delta^n_j-
    \delta^m_j\delta^n_i\,. 
\end{equation}
The commutator relations of $\mathfrak{so}(6)$ algebra read: 
\begin{equation}
    [L_{(m,n)}, L_{(p,q)}] = 
    \delta_{np}L_{(m,q)} + \delta_{nq}L_{(p,m)} + \delta_{mp}L_{(q,n)} + \delta_{mq}L_{(n,p)}\,.
\end{equation}

To prove the isomorphism between $SU(4)$ and the spin group 
$\mathrm{Spin}(6)$\, one usually proceeds from the isomorphism of their Lie algebras.  
Taking into account that $SO(6)$ has two 4-dimensional spinor representations, consider the following 1-to-1 mapping between the Fano basis of $\mathfrak{su}(4)$ and the Lie algebra basis 
(\ref{eq:basisso(6)}) of $\mathfrak{so}(6)$:
\footnote{Similar relations using the Gell-Mann basis for 
$\mathfrak{su}(4)$ algebra can be found in 
\cite{BarnesHamilton-CharltonLawrence2001}. 
} 
\begin{eqnarray}
  && \lambda_1 
  \leftrightarrow 
  L_{(1,5)}\,,\quad   \lambda_2 
  \leftrightarrow 
  L_{(1,3)}\,,\quad
   \lambda_3 
   \leftrightarrow 
   L_{(3,5)}\,,\quad
 \lambda_4 
 \leftrightarrow 
 L_{(2,4)}\,,\quad \lambda_5 
 \leftrightarrow 
 L_{(4,6)}\,,\\
&&   \lambda_6 
   \leftrightarrow 
   L_{(2,6)}\,,\quad
 \lambda_7
 \leftrightarrow 
 L_{(3,6)}\,,\quad   \lambda_8 
 \leftrightarrow 
 L_{(3,2)}\,,\quad
   \lambda_9
   \leftrightarrow 
   L_{(4,3)}\,,\quad
\lambda_{10}
\leftrightarrow 
L_{(6,5)}\,,\\  
&&
\lambda_{11} 
\leftrightarrow 
L_{(2,5)}\,,\quad
   \lambda_{12} 
   \leftrightarrow 
   L_{(5,4)}\,,\quad
 \lambda_{13}
 \leftrightarrow 
 L_{(6,1)}\,,\quad   \lambda_{14} \leftrightarrow L_{(2,1)}\,,\quad
   \lambda_{15} 
   \leftrightarrow 
   L_{(1,4)}\,. 
\end{eqnarray}
Under this mapping, the algebra $\mathfrak{so}(6)$ is decomposed into the following direct sum: 
\begin{equation}
\label{eq:so6dec}
    \mathfrak{so}(4) \oplus \mathfrak{t}^{''} \oplus \mathfrak{t}^{'} \oplus \mathfrak{t}\,,  
\end{equation} 
where 
$\mathfrak{t}, \mathfrak{t}^\prime$, and $\mathfrak{t}^{''}$ are three copies of the 3-dimensional  maximal Abelian  subalgebras of $\mathfrak{so}(6):$
\begin{equation}
\mathfrak{t}=\mbox{span}_{\mathbb{R}}\left(L_{(3,5)}, L_{(2,6)}, L_{(1,4)}\right)\,,
\qquad
\mathfrak{t}^{'}=\mbox{span}_{\mathbb{R}}\left(L_{(4,3)}, L_{(2,5)}, L_{(6,1)}\right)\,, 
\qquad
\mathfrak{t}^{''}=\mbox{span}_{\mathbb{R}}\left(L_{(1,5)}, L_{(2,4)}, L_{(3,6)}\right)\,.\nonumber 
\end{equation}
One can verify that under this mapping the commutator relations between the elements of decomposition (\ref{eq:so6dec}) reproduce the corresponding relations of the $\mathfrak{su(4)}$
algebra.
As a result of the exponential mapping, we arrive at 
the representation for $g \in \mathrm{SO}(6):$  
\begin{equation}
\label{eq:ort6}
g = \mathcal{R}_{{}_{\mathrm{SO(4)}}}
\, C_{\pi^{''}} 
O(\boldsymbol{\alpha})
C_{\pi^{'}} 
O(\boldsymbol{\beta})
C_{\pi}\,
O(\boldsymbol{\theta}) \,,
\end{equation} 
where $C_{\pi}, C_{\pi^{'}}, C_{\pi^{''}}$
denote the permutation matrices corresponding to the transpositions  
\(
\pi^{''} = (1,6,3,2,4,5)\,, \quad 
    \pi^{'} = (1,4,3,5,6)\,, \) and  
\(   \pi = (1,4,5,3,6)\,
\), respectively. 
% $R_{(1,6,3,2,4,5)} C_{(1,5,3)\times(2,4)} = C_{(1,4,3,5,6)} = R_{(1,6,5,3,4)}$, $R_{(1,5,3)\times(2,4)} C_{(1,2,4,3,6,5)} = C_{(1,4,5,3,6)} = R_{(1,6,3,5,4)}$\,. 
In (\ref{eq:ort6}) the 
first factor $\mathcal{R}_{{}_{\mathrm{SO(4)}}} \in \mathrm{SO}(4) \subset SO(6)\,$ 
and block-diagonal orthogonal matrices 
\begin{equation}
O(\boldsymbol{\chi}) = 
\begin{Vmatrix}
   R(\chi_1) & 0 &  0\\
    0& R(\chi_2)& 0\\
    0 &0&  R(\chi_3)
\end{Vmatrix} 
\end{equation}
describe rotations in three pairs of orthogonal  2-dimensional planes on angles $\chi_1, \chi_2, \chi_3\,$,   respectively,
\[
R(\chi_i)=
\begin{Vmatrix}
   \cos(\chi_i) & -\sin(\chi_i)\\
    \sin(\chi_i) &  \cos(\chi_i)
\end{Vmatrix}\,,
\quad i=1,2,3\,.
\]
The triplets $\boldsymbol{\alpha}, \boldsymbol{\beta}$, and 
$\boldsymbol{\theta}$ are definite over the direct product of 3 copies of the 3-cube, 
$
\boldsymbol{\alpha}\in 
[-\pi, \pi]^3\,,\
\boldsymbol{\beta}\in 
[-\pi, \pi]^3\,,
\
\boldsymbol{\theta} \in 
[-\pi, \pi]^3$.
Therefore, we ascertained that angles $\boldsymbol{\alpha}, \boldsymbol{\beta}$, and 
$\boldsymbol{\theta}$ from the three octahedrons represent the unitary covering counterparts of the rotation angles in the orthogonal matrix (\ref{eq:ort6}), in accord with the fact that $SU(4)$ is a double cover of $SO(6)$.

%%%%%%%%%%%%%
\section{Concluding remarks}
%%%%%%%%%%%%%%%%%%%%%%

In the paper we introduce a set of real coordinates on the $SU(4)$ group and determine their domain of definition allowing to describe almost all unitary transformations.  These coordinates emerge  as a result of factorization of the
special unitary transformation into left and right factors, corresponding to the left and right action of its subgroups,  $SU(2)\times SU(2)$ and maximal Abelian torus $T^3$, respectively. Formally,  it  was shown that an element of $SU(4)$ can be written as 
\begin{equation}
\label{eq:K3D}
g = k(\boldsymbol{u},\boldsymbol{v})\, H\,D(\boldsymbol{\alpha})\, H^{'}\,D(\boldsymbol{\beta})\,
H^{''}D(\boldsymbol{\theta})\,,
\end{equation}
where  $k$  belongs to the subgroup $K =SU(2)\times SU(2),$
represented as (\ref{eq:K2});
$D$ is the  diagonal matrix function from
(\ref{eq:diag}); while  $H, H^{'}$ and $ H^{''}$
are the following  Hadamard  matrices: 
\begin{equation}
%\scriptsize
    \label{eq:3Had}
    H=\frac{1}{2}\,
    \left(
\begin{array}{cccc}
 1 & 1 & -1 & -1 \\
 1 & -1 & -1 & 1 \\
 1 & -1 & 1 & -1 \\
 1 & 1 & 1 & 1 \\
\end{array}
\right),
\quad
H^{'}=\frac{1}{2}\,
    \left(
\begin{array}{cccc}
 -1 & 1 & 1 & 1 \\
 1 & 1 & -1 & 1 \\
 1 & -1 & 1 & 1 \\
 1 & 1 & 1 & -1 \\
\end{array}
\right),
\quad
H^{''}=\frac{1}{2} \, \left(
\begin{array}{cccc}
 -1 & -1 & -1 & 1 \\
 1 & 1 & -1 & 1 \\
 -1 & 1 & 1 & 1 \\
 1 & -1 & 1 & 1 \\
\end{array}
\right)\,.
\end{equation} 

To our best knowledge, the suggested decomposition (\ref{eq:K3D}) is a new one and it is well adapted to  the  computational issues  where SVD is applicable. 
Particularly, in the forthcoming publications we plan to exploit the proposed parameterization for studies of  the entanglement characteristics of 2-qubits. 
According to (\ref{eq:K3D}), the point in double coset $\mathcal{A}=SU(2)\times SU(2)\backslash SU(4)/T^3 $ can be parameterized as 
\begin{equation}
\mathcal{A}= 
H\,D(\boldsymbol{\alpha}) H^{'}\, D(\boldsymbol{\beta})\,
H^{''} \,.   
\end{equation}
Hence,  using parameterization (\ref {eq:K3D}) for the diagonalizing unitary matrix of a 2-qubit density matrix $\varrho= g\,\mbox{diag}(r_1, r_2, r_3, r_4 )\,g^\dagger$, we identify  triples $\boldsymbol{\alpha}$, $\boldsymbol{\beta}$ and eigenvalues $r_1, r_2, r_3, r_4$ with  nine independent $SU(2)\times SU(2)$\--invariants characterizing  the  entanglement space of a qubit pair.   
%%%%%%%%%%%%%%%%%%%
\section*{Acknowledgments}
%%%%%%%%%%%%%%%%%%%%%%%%%%%%%

The research of A.K and D.M. was supported in part by the "Bulgaria-JINR" collaborative grant, the work of A.T. was partially funded by the Higher Education and Science Committee of MESCS RA (research project № 23/2IRF-1C003).

%%%%%%%%%%%%%%%%%%%%%%%%%%%%%%%%%%%%%%

%%%%%%%%%%%%%%%%%%%%%%% APPENDIX %%%%%%%%%%%%%%%%%%
\section{Appendix}
\label{Appendix}
\appendix
%%%%%%%%%%%%%%%%%%%%%%%%%%%%%%%%%

\paragraph{Fano basis of the Lie algebra $\mathfrak{su}(4)$.}
Starting from  the fundamental representation of the $\mathfrak{su}(2)$ algebra  with basis $\{ \frac{1}{2\imath}\sigma_1, \frac{1}{2\imath}\sigma_2, \frac{1}{2\imath}\sigma_3\}$ 
in terms of the Pauli sigma-matrices, let us build  tensor products   
\[
\sigma_{i0}= \frac{1}{2\imath}\,
\sigma_i\otimes \mathbb{I}_2\,,\quad
\sigma_{0i}=\frac{1}{2\imath}\,
\mathbb{I}_2\otimes\sigma_i\,,
\quad 
\sigma_{ij}= \frac{1}{2\imath}\,
\sigma_i\otimes \sigma_j\,, \qquad 
i,j =1,2,3\,, 
\]
and set the ordered 15-tuple
$\boldsymbol{\lambda}=\{\lambda_1,\lambda_2,\dots,\lambda_{15}\}\,$
as follows: 
\begin{equation} 
\boldsymbol{\lambda}=\{
 \sigma_{10}, \sigma_{20}, \sigma_{30},  \sigma_{01}, \sigma_{02}, \sigma_{03},
  \sigma_{11}, \sigma_{12}, \sigma_{13},
   \sigma_{21}, \sigma_{22}, \sigma_{23}, 
    \sigma_{31}, \sigma_{32}, \sigma_{33}
\}\,.
\end{equation} 
Matrices  $\boldsymbol{\lambda}$ comprise the anti-Hermitian Fano basis of the Lie algebra $\mathfrak{su}(4)\,,$
which is orthonormal in the sense that 
$\mbox{Tr}(\lambda_i\lambda_j )=-1\,.$
The explicit form of the Fano basis matrices is:    
%%%

{\scriptsize
\bigskip
\noindent$
\lambda_1\,=\displaystyle{\frac{1}{2\imath}}
\begin{Vmatrix}
        0 & 0 & 1 & 0 \\
        0 & 0 & 0 & 1 \\
        1 & 0 & 0 & 0 \\
        0 & 1 & 0 & 0 \\
 \end{Vmatrix},
\qquad
\lambda_2=\displaystyle{\frac{1}{2\imath}}
 \begin{Vmatrix}
      0 & 0 & -i & 0 \\
      0 & 0 & 0 & -i \\
      i & 0 & 0 & 0 \\
      0 & i & 0 & 0 \\
    \end{Vmatrix},
   \qquad
  \lambda_3=\displaystyle{\frac{1}{2\imath}}
  \begin{Vmatrix}
      1 & 0 & 0 & 0 \\
      0 & 1 & 0 & 0 \\
      0 & 0 & -1 & 0 \\
      0 & 0 & 0 & -1 \\
     \end{Vmatrix}, 
    \qquad 
     \lambda_4\,=\displaystyle{\frac{1}{2\imath}}
\begin{Vmatrix}
       0 & 1 & 0 & 0 \\
       1 & 0 & 0 & 0 \\
       0 & 0 & 0 & 1 \\
       0 & 0 & 1 & 0 \\
      \end{Vmatrix},
\qquad
\lambda_5=\displaystyle{\frac{1}{2\imath}}
 \begin{Vmatrix}
      0 & -i & 0 & 0 \\
      i & 0 & 0 & 0 \\
      0 & 0 & 0 & -i \\
      0 & 0 & i & 0 \\
     \end{Vmatrix},
$
\bigskip
\bigskip

\noindent$ 
\lambda_6=\displaystyle{\frac{1}{2\imath}}
      \begin{Vmatrix}
      1 & 0 & 0 & 0 \\
      0 & -1 & 0 & 0 \\
      0 & 0 & 1 & 0 \\
      0 & 0 & 0 & -1 \\
     \end{Vmatrix},
\qquad
     \lambda_7\,=\displaystyle{\frac{1}{2\imath}} 
     \begin{Vmatrix}
      0 & 0 & 0 & 1 \\
      0 & 0 & 1 & 0 \\
      0 & 1 & 0 & 0 \\
      1 & 0 & 0 & 0 \\
     \end{Vmatrix},
\qquad
  \lambda_8=\displaystyle{\frac{1}{2\imath}} 
  \begin{Vmatrix}
      0 & 0 & 0 & -i \\
      0 & 0 & i & 0 \\
      0 & -i & 0 & 0 \\
      i & 0 & 0 & 0 \\
     \end{Vmatrix},
\qquad
   \lambda_9=\frac{i}{2}  \begin{Vmatrix}
      0 & 0 & 1 & 0 \\
      0 & 0 & 0 & -1 \\
      1 & 0 & 0 & 0 \\
      0 & -1 & 0 & 0 \\
     \end{Vmatrix}, 
\quad
\lambda_{10}=\displaystyle{\frac{1}{2\imath}}
\begin{Vmatrix}
      0 & 0 & 0 & -i \\
      0 & 0 & -i & 0 \\
      0 & i & 0 & 0 \\
      i & 0 & 0 & 0 \\
     \end{Vmatrix},
$
\bigskip
\bigskip

\noindent$  
\lambda_{11}=\displaystyle{\frac{1}{2\imath}}
\begin{Vmatrix}
      0 & 0 & 0 & -1 \\
      0 & 0 & 1 & 0 \\
      0 & 1 & 0 & 0 \\
      -1 & 0 & 0 & 0 \\
     \end{Vmatrix},
     \quad
\lambda_{12}=\displaystyle{\frac{1}{2\imath}}
 \begin{Vmatrix}
      0 & 0 & -i & 0 \\
      0 & 0 & 0 & i \\
      i & 0 & 0 & 0 \\
      0 & -i & 0 & 0 \\
     \end{Vmatrix},   
\quad
\lambda_{13}=\displaystyle{\frac{1}{2\imath}}
\begin{Vmatrix}
      0 & 1 & 0 & 0 \\
      1 & 0 & 0 & 0 \\
      0 & 0 & 0 & -1 \\
      0 & 0 & -1 & 0 \\
     \end{Vmatrix},
\quad
\lambda_{14}=
\displaystyle{\frac{1}{2\imath}}
\begin{Vmatrix}
      0 & -i & 0 & 0 \\
      i & 0 & 0 & 0 \\
      0 & 0 & 0 & i \\
      0 & 0 & -i & 0 \\
     \end{Vmatrix},
     \quad
\lambda_{15}=
\displaystyle{\frac{1}{2\imath}}
\begin{Vmatrix}
      1 & 0 & 0 & 0 \\
      0 & -1 & 0 & 0 \\
      0 & 0 & -1 & 0 \\
      0 & 0 & 0 & 1 \\
\end{Vmatrix}. 
$
}

In the Fano basis, the skew-symmetric constants
$\mathrm{f}^k_{\,ij}$  
determining the commutator relations  
\begin{equation}
[\lambda_i\,, \lambda_j ]= \sum_{k=1}^{15}\,\mathrm{f}^{k}_{\,ij}\,\lambda_k\,,
\end{equation}
take the values $0$ and $\pm 1\,.$ All pairwise commutators between basis elements can be extracted from Table \ref{T:CommRelsu(4)}. 

The product of any two Fano matrices can be expanded over the basis elements supplemented by a unit $4\times 4$ matrix $\mathbb{I}_4$,
\begin{eqnarray}
&\lambda_{i} \lambda_{j}=-\frac{1}{4}\,\delta_{i j}\,\mathbb{I}_4 + \frac{1}{2}\,\left(\mathrm{f}^k_{\, i j}+
\imath\,\mathrm{d}^k_{\, i j}\right)\lambda_{k}\,.
\end{eqnarray}
Non-vanishing  symmetric coefficients   
$\mathrm{d}^k_{\, i j}= - \imath\, \mbox{Tr}(\{\lambda_{i}\,, \lambda_{j}\} \lambda_{k})$ are listed (up to permutations) in Table \ref{T:SymCoeffSU(4)}.
%%%%%%%%%%%%%%%%%%%%%%%%%
%%%  Table 1. Algebra su(4) %%%%%%%%%%%%%%
{
%\begin{landscape}
\begin{table}[h]
\scriptsize
\centering
%\rowcolors{5-8}{yellow}{pink}
%\rowcolors{2}{pink}{pink}
\begin{tabular}{|>{\columncolor{green}} c| c |c|>{\columncolor{white}} c||c|c|c||c|c|c||c|c|c||c|c
|>{\columncolor{white}} c |c}
\hline
\cellcolor{white}&\multicolumn{3}{|c||}{\cellcolor{green}$\boldsymbol{\mathfrak{a}}$}& 
\multicolumn{3}{c||}{\cellcolor{yellow}$\boldsymbol{\mathfrak{a}^\prime}$}&
\multicolumn{6}{c||}{\cellcolor{pink}$\boldsymbol{\mathfrak{k}}$} &
\multicolumn{3}{c|}{\cellcolor{lightgray}$\mathfrak{t}$}\\
\hline
\hline
\rowcolor{white}& 
$\cellcolor{green}\lambda_1 $&
$\cellcolor{green}\lambda_4$&
\cellcolor{green}$\lambda_7$&
%%%%%
$\cellcolor{yellow}\lambda_9$ &
\cellcolor{yellow}$\lambda_{11}$&
$\cellcolor{yellow}\lambda_{13}$ &
%%%%
$\cellcolor{pink}\lambda_2$ &
$\cellcolor{pink}\lambda_8$ &
$\cellcolor{pink}\lambda_{14}$ &
$\cellcolor{pink}\lambda_{5}$ &
$\cellcolor{pink}\lambda_{10}$ &
$\cellcolor{pink}\lambda_{12}$&
%%%%
\cellcolor{lightgray}$\lambda_3 $&
\cellcolor{lightgray}$\lambda_6 $&
\cellcolor{lightgray}$\lambda_{15}$ \\ 
\hline\hline
%1
$\cellcolor{green}\lambda_1 $& 
0  & 
0 & 
0 & 
0 &
-$\cellcolor{pink}\lambda_{14}$&
$\cellcolor{pink}\lambda_{10}$&
-$\cellcolor{lightgray}\lambda_{3}$&
0&
\cellcolor{yellow}$\lambda_{11}$&
0&
-$\cellcolor{yellow}\lambda_{13}$&
$\cellcolor{lightgray}\lambda_{15}$&
$\cellcolor{pink}\lambda_{2}$&
0&
$\cellcolor{pink}\lambda_{12}$
\\ \hline
%4
$\cellcolor{green}\lambda_4$ & 
0 & 
0 &
0&
$\cellcolor{pink}\lambda_8$&
-$\cellcolor{pink}\lambda_{12}$&
0&
0&
$-\cellcolor{yellow}\lambda_9$&
$-\cellcolor{lightgray}\lambda_{15}$&
-\cellcolor{lightgray}$\lambda_{6}$&
0&
$\cellcolor{yellow}\lambda_{11}$&
0&
$\cellcolor{pink}\lambda_{5}$&\cellcolor{pink}$\lambda_{14}$
\\\hline
% 7
%\rowcolor{white}
$\cellcolor{green}\lambda_7$&
0 & 
0&
0& 
\cellcolor{pink} $\lambda_5$&
0&
$\cellcolor{pink}\lambda_2$ &
-$\cellcolor{yellow}\lambda_{13}$&
-\cellcolor{lightgray} $\lambda_6$&
0&
$-\cellcolor{yellow}\lambda_9$&
$-\cellcolor{lightgray}\lambda_3$&
0&
$\cellcolor{pink}\lambda_8$&
0&0\\
\hline
\hline
%%%%%%%%%%%%%%%%%%%%%%%%%%%%%%%%%%%%%
% 9

$\cellcolor{yellow}\lambda_9$&
0&
-$\cellcolor{pink}\lambda_{8}$&
-\cellcolor{pink}$\lambda_{5}$&
0&
0&
0&
-$\cellcolor{lightgray}\lambda_{15}$&
$\cellcolor{green}\lambda_4$&
0&
\cellcolor{green}$\lambda_{7}$&
0&
-$\cellcolor{lightgray}\lambda_3$&
$\cellcolor{pink}\lambda_{12}$&
0&
\cellcolor{pink}$\lambda_2$\\
\hline
% 11
%\rowcolor{white}
$\cellcolor{yellow}\lambda_{11}$ &\cellcolor{pink}$\lambda_{14}$&
\cellcolor{pink}$\lambda_{12}$&
0&
0&
0&
0&
0&
\cellcolor{lightgray}$\lambda_3$&
\cellcolor{green}$\lambda_1$&
0&
\cellcolor{lightgray}$\lambda_6$&
\cellcolor{green}-$\lambda_4$&
$-\cellcolor{pink}\lambda_8$&0&0\\
\hline
%13
$\cellcolor{yellow}\lambda_{13} $&-$\cellcolor{pink}\lambda_{10}$&
0&
$-\cellcolor{pink}\lambda_2$&
0&
0&
0&
$\cellcolor{green}\lambda_7$&
0&
-\cellcolor{lightgray}$\lambda_{6}$&
$-\cellcolor{lightgray}\lambda_{15}$&
\cellcolor{green}$\lambda_{1}$&
0&
0&
$\cellcolor{pink}\lambda_{14}$&\cellcolor{pink}$\lambda_5$\\
\hline
\hline
\hline
%%%%%%%%%%%%%%%%%%%%%%%%%%%%%
%2 
%\rowcolor{pink}
$\cellcolor{pink}\lambda_2 $ & $\cellcolor{lightgray}\lambda_3$ &0&
\cellcolor{yellow}-$\lambda_{13}$&
\cellcolor{lightgray}$\lambda_{15}$&
0&
$-\cellcolor{green}\lambda_{7}$&
0&
$\cellcolor{pink}\lambda_{14}$&
-$\cellcolor{pink}\lambda_{8}$&
0&
0&
0&
-$\cellcolor{green}\lambda_{1}$&
0&
\cellcolor{yellow}$-\lambda_{9}$\\
\hline
% 8 
%\rowcolor{white}
$\cellcolor{pink}\lambda_8$&
0&
$\cellcolor{yellow}\lambda_{9}$&
$\cellcolor{lightgray}\lambda_{6}$ &
-\cellcolor{green}$\lambda_4$&
$-\cellcolor{lightgray}\lambda_3$&
0&
$-\cellcolor{pink}\lambda_{14}$&
0&
\cellcolor{pink}$\lambda_2$&
0&
0&
0&
\cellcolor{yellow}$\lambda_{11}$
&-\cellcolor{green}$\lambda_7$& 0\\
\hline
%14
$\cellcolor{pink}\lambda_{14}$ &-$\cellcolor{yellow}\lambda_{11}$&
$\cellcolor{lightgray}\lambda_{15}$&
0&
0&
$\cellcolor{green}\lambda_1$
&-$\cellcolor{lightgray}\lambda_{6}$&
\cellcolor{pink}$\lambda_8$&
-$\cellcolor{pink}\lambda_2$&
0&
0&
0&
0&
0&
-$\cellcolor{yellow}\lambda_{13}$&
\cellcolor{green}-$\lambda_4$\\
\hline
\hline
% 5
 $\cellcolor{pink}\lambda_5$&0 &\cellcolor{lightgray}$\lambda_6$&
 $\cellcolor{yellow}\lambda_9$&
 $-\cellcolor{green}\lambda_7$&
 0&
 $
\cellcolor{lightgray}\lambda_{15}$&
0&
0&
0&
0&
$\cellcolor{pink}\lambda_{12}$&
-$\cellcolor{pink}\lambda_{10}$&
0&
-\cellcolor{green}$\lambda_4$&
\cellcolor{yellow}$-\lambda_{13}$\\
\hline 
% 10
%\rowcolor{white}
$\cellcolor{pink}\lambda_{10}$&
\cellcolor{yellow}$\lambda_{13}$&
0&
$\cellcolor{lightgray}\lambda_3$&
0&
\cellcolor{lightgray}-$\lambda_{6}$&$-\cellcolor{green}\lambda_{1}$&
0&
0&
0&
-\cellcolor{pink}$\lambda_{12}$&
0&
\cellcolor{pink}$\lambda_5$&
\cellcolor{green}-$\lambda_7$&
$\cellcolor{yellow}\lambda_{11}$
&0\\\hline
% 12
$\cellcolor{pink}\lambda_{12}$&$\cellcolor{lightgray}\lambda_{15}$&-$\cellcolor{yellow}\lambda_{11}$  &
0&
$\cellcolor{lightgray}\lambda_{3}$&
$\cellcolor{green}\lambda_{4}$&0&
0&
0&
0&
\cellcolor{pink}$\lambda_{10}$&
-$\cellcolor{pink}\lambda_5$&
0&
$-\cellcolor{yellow}\lambda_9$&
0&
\cellcolor{green}-$\lambda_1$\\
\hline
\hline
\hline
%%%%%%%%%%%%%%%%%%%%%%%%%%%
%3
%\rowcolor{white}
$\cellcolor{lightgray}\lambda_3$&
-\cellcolor{pink}$\lambda_2$&
0&
-\cellcolor{pink}$\lambda_{10}$&
-\cellcolor{pink}$\lambda_{12}$&
\cellcolor{pink}$\lambda_{8}$&
0&
$\cellcolor{green}\lambda_{1}$&
-\cellcolor{yellow}$\lambda_{11}$&
0&
0&
\cellcolor{green}$\lambda_7$&
\cellcolor{yellow}$\lambda_9$&
0&
0&
0\\
\hline
% 6
%\rowcolor{white}
$\cellcolor{lightgray}\lambda_6$&
0&
-\cellcolor{pink} $\lambda_5$&
\cellcolor{pink}-$\lambda_8$&
0&
\cellcolor{pink}$\lambda_{10}$&
-$\cellcolor{pink}\lambda_{14}$&
0&
$\cellcolor{green}\lambda_{7}$&
\cellcolor{yellow}$\lambda_{13}$&
$\cellcolor{green}\lambda_{4}$&
-\cellcolor{yellow}$\lambda_{11}$&
0&
0&
0&
0\\
\hline
%15
%\rowcolor{white}
$\cellcolor{lightgray}\lambda_{15}$ & \cellcolor{pink}-$\lambda_{12} $ &
 -\cellcolor{pink}$\lambda_{14}$ 
 &0 &
 -\cellcolor{pink}-$\lambda_{2}$ &
 0 &
 -$\cellcolor{pink}\lambda_5$&
 $\cellcolor{yellow}\lambda_{9}$&
 0 &
$\cellcolor{green}\lambda_{4}$&
 $\cellcolor{yellow}\lambda_{13}$&
0&
\cellcolor{green}$\lambda_1$&
0 &
0&
0\\
\hline
\end{tabular}
\bigskip
\bigskip
\caption{Commutator relations between  elements of the Fano basis $\boldsymbol{\lambda}$ collected in accordance with  the direct sum  decomposition of the $\mathfrak{su(4)}$ algebra;   
$\mathfrak{k}=\mathfrak{su(2)}\oplus\mathfrak{su(2)}$ ({\color{pink}pink}), the maximal torus $\mathfrak{t}$ {\color{lightgray}(lightgray)}, and two Abelian subalgebras $\mathfrak{a}$ {\color{yellow}(yellow)} and $\mathfrak{a}^\prime$ {\color{green} (green)}.}
\label{T:CommRelsu(4)}
\end{table}
%\end{landscape}
}
%%%%%%%%%%% Symmetric coefficients
{
\begin{table}[ht]
\scriptsize
\begin{center}
\begin{tabular}{|| c || c ||}
\hline
\hline\begin{tabular}[t]{c c c}&\cellcolor{lightgray} $d_{i j k}=1$&\\
 \hline
 \hline
 i & j & k\\
 \hline
 1 & 4 & 7 \\
 1 & 5 & 8 \\
 1 & 6 & 9 \\
 2 & 4 & 10 \\
 2 & 5 & 11 \\
 2 & 6 & 12 \\
 3 & 4 & 13 \\
 3 & 5 & 14 \\
 3 & 6 & 15 \\
 7 & 12 & 14 \\
 8 & 10 & 15 \\
 9 & 11 & 13 \\
\end{tabular}
& 
\begin{tabular}[t]{c c c}
 &\cellcolor{lightgray}$d_{i j k}=-1$&\\
 \hline
 \hline
 i & j & k\\
 \hline
 7 & 11 & 15 \\
 8 & 12 & 13 \\
 9 & 10 & 14 \\
\end{tabular}\\
\hline
\hline
\end{tabular}
\bigskip
\bigskip
\caption{The nonzero symmetric structure constants of the $\mathfrak{su(4)}$ algebra up to cyclic permutations.}
\label{T:SymCoeffSU(4)}
\end{center}
\end{table}
}

\newpage

%%%%%%%%%%%%%%%%%%%%%%%
%%%%%% REFERENCES 
%%%%%%%%%%%%%%%%%%%%%%

%

\begin{thebibliography}{99}
%
\bibitem{Helgason1978}
Helgason, S.: 
Differential Geometry, Lie Groups, and Symmetric Spaces. 
Academic Press (1978)
%
\bibitem{Gilmore2012}
Gilmore, R.: 
Lie Groups, Lie Algebras, and Some of their Applications. Courier Corporation (2012)
%
\bibitem{Murnaghan1952}
Murnaghan, F.D.: 
On convenient system of parameters for the unitary group. 
The Proceedings of the National Academy of Sciences 
\textbf{38}(2), 127–129 (1952)
%
\bibitem{Murnghan1962}
Murnaghan, F.D.: 
The unitary and rotation groups.
Spartan, Washington (1962)
%
\bibitem{Givens1958}
Givens, W.: 
Computation of Plain Unitary Rotations Transforming a General Matrix to Triangular Form. 
Journal of the Society for Industrial and Applied Mathematics \textbf{6}(1), 26–50 (1958)
%
\bibitem{GolubVanLoan2012}
Golub, G.H, Van Loan, C.F.: 
Matrix computations. JHU Press, 4th edition, (2012)
%
\bibitem{Wigner1968}
Wigner, E.P.: 
On a generalization of Euler’s angles. 
In: Group Theory and Its Applications. 
Elsevier, pp. 119–129 (1968)
%
\bibitem{Macfarlane1968}
Macfarlane, A.J.: 
Description of the symmetry group $SU(3)/Z_3$.
Commun. Math. Phys. \textbf{11}, 91-98 (1968)
%
\bibitem{MichelRadicati1973}
Michel, L., Radicati, L.A. 
The geometry of the octet. 
Ann. Inst. Henri Poincare, \textbf{XVIII}(3), 185-214
(1973) 
%
\bibitem{Bincer1990}
Bincer, A.M.: 
Parameterization of $SU(n)$ with $n-1$ orthonormal vectors.  J. Math. Phys. \textbf{31}, 563-567 (1990) 
%
\bibitem{Dita2003}
Dita, P.: 
Factorization of unitary matrices. 
J. Phys. A: Math. Gen. \textbf{36}(11), 2781 (2003) 
%		
\bibitem{TilmaSudarshan} 
Tilma, T., Sudarshan, E.C.G.: 
Generalized Euler angle parametrization for $U(N)$ with application to $SU(N)$ coset volume measures.
Journal of Geometry and Physics, \textbf{52} (3), 263-283 (2004)
%
\bibitem{BertiniCacciatoriCerchiai} 
Bertini, S., Cacciatori, S.L, B. L. Cerchiai, B.L.:
On the Euler angles for $SU(N)$. 
J. Math. Phys. \textbf{47}(4), 043510 (2006)
%		
\bibitem{CacciatoriDallaPiazzaScotti2017} 
Cacciatori, S.L., Dalla Piazza, F., Scotti, A.: 
Compact Lie groups: Euler constructions and generalized Dyson conjecture. 
Trans. Amer. Math. Soc. \textbf{369}, 4709-4724 (2017)
%
\bibitem{SpenglerHuberHiesmayr}
Spengler, C., Huber M., Hiesmayr B.: 
A composite parameterization of unitary groups, density matrices and subspaces.
J. Phys. A: Math. Theor. \textbf{43}, 385306 (2010) 
%
\bibitem{Charzynski2005}
Charzynski, S., Kijowski, J., Rudolph, G., Schmidt, M.:
On the stratified classical configuration space of lattice QCD.
J. Geom. Phys. \textbf{55}, 137-178 (2005)
%
\bibitem{EdelmanJeong2022}
Edelman, A., Jeong, S.: 
On the Cartan decomposition for classical random matrix ensembles. 
J. Math. Phys. \textbf{63}, 061705 (2022)
%
\bibitem{Matsuki1995}
Matsuki, T.: 
Double coset decompositions of algebraic groups arising from two involutions. 
Journal of Algebra \textbf{175}(3), 865-925 (1995)
%
\bibitem{Matsuki1997}
Matsuki T.: 
Double coset decompositions of reductive Lie groups arising from two involutions,
Journal of algebra \textbf{197}, 49-91 (1997)
%
\bibitem{Miebach2007}
Miebach, C.:  
Matsuki's double coset decomposition via gradient maps. Journal of Lie Theory \textbf{18}(3), 555-580 (2008);
arXiv:0712.3384[math.RT] (2007)
%
\bibitem{Kobayashi2007}
Kobayashi, T.:  
A generalized Cartan decomposition for the double coset space $(U (n_1)\times U (n_2) \times U (n_3)) 
\backslash U (n)/(U (p) \times U (q))$. 
Journal of the Mathematical Society of Japan \textbf{59}(3), 
669-691, (2007)
%
\bibitem{EdelmanJeong2023} 
Edelman, A., Jeong, S.:  
Fifty three matrix factorizations: A systematic approach. 
SIAM Journal on Matrix Analysis and Applications 
\textbf{44}(2), 415-480 (2023)
%
\bibitem{Postnikov1986}
Postnikov, M.M.: 
Lectures in Geometry: Lie Groups and Lie Algebras. Semester 5. English translation, Mir Publishers (1986) 
%
\bibitem{Lee2013}
Lee, J.M.: 
Introduction to smooth manifolds. 
Graduate Texts in Mathematics, Second Edition, Springer, (2013)
%
\bibitem{FujiiSuzuki2007}
Fujii, K., Suzuki, T.: 
On the magic matrix by Makhlin and the BCH formula for SO(4). arXiv:quant-ph/0610009 (2007)
%
\bibitem{Higham2008}
Higham, N.J.: 
Functions of Matrices. 
Society for Industrial and Applied Mathematics (SIAM), Philadelphia, PA (2008) 
%
\bibitem{BarnesHamilton-CharltonLawrence2001}
 Barnes, K.J., Hamilton-Charlton, K.J., Lawrence, T.R.: 
 How orbits of $SU(N)$ can describe rotations in $SO(6)$.
J. Phys. A: Math. Gen. \textbf{34}, 10881 (2001)
%
\end{thebibliography}
\end{document}